\documentclass{amsart}
\usepackage{graphicx} % Required for inserting images

\usepackage{hyperref}   % www references
\usepackage[square,numbers]{natbib}     % bibliography
\usepackage{amsthm}     % ams theorem environments
\usepackage{amsmath}    % ams math environments
\usepackage{amssymb}    % ams symbols
\usepackage{amsxtra}    % ams extra symbols
\usepackage{eucal}      % eu calligraphic fonts
\usepackage[bbgreekl]{mathbbol}   % black board symbols
\usepackage{todonotes}
\usepackage{comment}
\usepackage{tikz}

\usepackage{diagbox}
\usepackage{color, colortbl}
\usepackage{array}
\usepackage{varwidth} %for the varwidth minipage environment

\definecolor{Gray}{gray}{0.9}
\definecolor{LightGray}{gray}{0.95}

\newcolumntype{g}{>{\columncolor{Gray}}c}
\newcolumntype{M}{V{3cm}}
\newcolumntype{D}{V{4cm}}
\newcolumntype{H}{V{4cm}}
\newcolumntype{F}{p{.45cm}} % fixed width

\usetikzlibrary{matrix}
\usetikzlibrary{patterns, shapes}
\usetikzlibrary{arrows}
\usetikzlibrary{calc,3d}
\usetikzlibrary{decorations,decorations.pathmorphing, decorations.pathreplacing}
\usetikzlibrary{through}
\tikzset{every picture/.style={baseline=-.65ex}}
\tikzset{ext/.style={circle, draw,inner sep=1pt},int/.style={circle,draw,fill,inner sep=1pt},nil/.style={inner sep=1pt}}
\tikzset{exte/.style={circle, draw,inner sep=3pt},inte/.style={circle,draw,fill,inner sep=3pt}}
\tikzset{diagram/.style={matrix of math nodes, row sep=3em, column sep=2.5em, text height=1.5ex, text depth=0.25ex}}
\tikzset{diagram2/.style={matrix of math nodes, row sep=0.5em, column sep=0.5em, text height=1.5ex, text depth=0.25ex}}
\tikzset{every loop/.style={draw}}

\newcommand{\tadpi}{
\begin{tikzpicture}[baseline=-.65ex]
\node[int] (v) at (0,0) {};
\draw (v) edge[loop] (v) edge +(0,-.5);
\end{tikzpicture}
}

\newcommand{\thetagr}{
\begin{tikzpicture}[baseline=-.65ex]
\node[int] (v) at (0.3,0) {};
\node[int] (w) at (-.3,0) {};
\draw (v) edge[bend left] (w) edge[bend right] (w) edge (w);
\end{tikzpicture}
}

\newcommand{\thetagi}{
\begin{tikzpicture}[baseline=-.65ex]
\node[int] (v) at (0.3,0) {};
\node[int] (w) at (-.3,0) {};
\draw (v) edge[bend left] (w) edge[bend right] (w) edge (w) edge +(.5,0);
\end{tikzpicture}
}

\newcommand{\thetagii}{
\begin{tikzpicture}[baseline=-.65ex]
\node[int] (v) at (0.3,0) {};
\node[int] (w) at (-.3,0) {};
\draw (v) edge[bend left] (w) edge[bend right] (w) edge +(.5,0) (w) edge +(-.5,0);
\end{tikzpicture}
}

\usepackage{tikz-cd}

\theoremstyle{plain}

\newtheorem{thm}[subsection]{Theorem}

\newtheorem{prop}[subsection]{Proposition}

\newtheorem{cor}[subsection]{Corollary}

\newcommand{\G}{\mathsf{G}}
\newcommand{\GC}{\mathsf{GC}}

\newcommand{\lp}{\text{-loop}}

\newcommand{\Com}{\mathrm{Com}}
\newcommand{\COp}{\mathsf{C}}
\newcommand{\POp}{\mathsf{P}}

\DeclareMathOperator{\Feyn}{Feyn}
\DeclareMathOperator{\Fe}{Fe}
\newcommand{\Q}{\mathbb{Q}}
\newcommand{\kk}{\mathfrak{k}}

\title{On the bridgeless graph complex}

\author{Thomas Willwacher}
\address{Department of Mathematics \\ ETH Zurich \\
R\"amistrasse 101 \\
8092 Zurich, Switzerland}
\email{thomas.willwacher@math.ethz.ch}

\date{}

\begin{document}

\begin{abstract}
    We discuss the cohomology of the bridgeless graph complex, that is, the subcomplex of the Kontsevich graph complex spanned by bridgeless graphs. 
\end{abstract}

\maketitle

\section{Introduction and main results}

Let $\G_n$ be the Kontsevich graph complex, spanned by $\geq3$-valent connected graphs, possibly with self-edges or multiple edges. It is equipped with a differential $d$ given by edge contraction.
We say that a connected graph is \emph{bridgeless} (or 2-edge-connected), if removing any edge does not split the graph into multiple connected components. For example, the first two of the following graphs are bridgeless, the third is not:
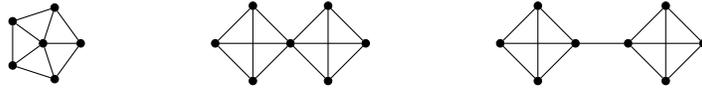
\begin{align*}
    &\begin{tikzpicture}
        \node[int] (v1) at (0:.5) {};
        \node[int] (v2) at (72:.5) {};
        \node[int] (v3) at (144:.5) {};
        \node[int] (v4) at (-144:.5) {};
        \node[int] (v5) at (-72:.5) {};
        \node[int] (c) at (0,0) {};
        \draw (c) edge (v1) edge (v2) edge (v3) edge (v4) edge (v5)
        (v1) edge (v2) edge (v5) (v3) edge (v2) edge (v4) (v4) edge (v5);
    \end{tikzpicture}
    &  &\begin{tikzpicture}
        \node[int] (v1) at (0,0) {};
        \node[int] (v2) at (1,0) {};
        \node[int] (v3) at (-1,0) {};
        \node[int] (v4) at (.5,.5) {};
        \node[int] (v5) at (.5,-.5) {};
        \node[int] (v6) at (-.5,.5) {};
        \node[int] (v7) at (-.5,-.5) {};
        \draw (v1) edge (v2) edge (v3) edge (v4) edge (v5) edge (v6) edge (v7)
        (v6) edge (v7) (v4) edge (v5)
        (v2) edge (v4) edge (v5)
        (v3) edge (v6) edge (v7);
    \end{tikzpicture}
    &  &\begin{tikzpicture}
        \node[int] (v1a) at (0,0) {};
        \node[int] (v1b) at (0.7,0) {};
        \node[int] (v2) at (1.7,0) {};
        \node[int] (v3) at (-1,0) {};
        \node[int] (v4) at (1.2,.5) {};
        \node[int] (v5) at (1.2,-.5) {};
        \node[int] (v6) at (-.5,.5) {};
        \node[int] (v7) at (-.5,-.5) {};
        \draw (v1b) edge (v2) edge (v1a) edge (v4) edge (v5)
        (v1a) edge (v3)  edge (v6) edge (v7)
        (v6) edge (v7) (v4) edge (v5)
        (v2) edge (v4) edge (v5)
        (v3) edge (v6) edge (v7);
    \end{tikzpicture}
    %\caption{\label{fig:graph examples} }
\end{align*}
Contracting an edge of a bridgeless graph results in a bridgeless graph.
Hence the subspace
\[
\G_n^{bl} \subset  \G_n
\]
spanned by bridgeless graphs is a subcomplex of the Kontsevich graph complex that we call the bridgeless graph complex. The author has been asked the following question on several occasions:

\medskip 

{\em What is the cohomology $H(\G_n^{bl})$ of the bridgeless graph complex?}

\medskip

This question has a fairly simple partial answer, in that there is a link to the cohomology of the (usual) Kontsevich graph complex with external legs. This result may also be obtained in significantly larger generality. To this end, let $\POp$ be a cyclic operad. Then we may define the Feynman transform $\Feyn(\POp)$ of $\POp$ as the $\kk$-modular cooperad whose cooperations in arity $r$ and genus $g$ are 
\begin{equation}\label{equ:feyn def}
\Feyn(\POp)(g,r)= \left(\bigoplus_{\Gamma} \bigotimes_{v\in V(\Gamma) } 
\POp(\mathrm{star}_v) \otimes \Q[1]^{\otimes E(\Gamma)}\right)/{\sim}
\end{equation}
Here the sum is over graphs of genus $g$ with $r$ external legs, $V(\Gamma)$ is the set of vertices of $\Gamma$, and one quotients by graph isomorphisms. 
Elements of $\Feyn(\POp)(g,r)$ can be understood as linear combinations of graphs of loop order $g$ with $r$ legs whose vertices are decorated by elements of $\POp$.
The modular operadic cocompositions are by cutting the graph, removing (at least) one edge.
\begin{equation}\label{equ:cocomp ex}
\begin{tikzpicture}[scale=.7]
                \node[int] (v1) at (-.5,0) {};
                \node[int] (v2) at (0,.5) {};
                \node[int] (v3) at (0,-.5) {};
                \node[int] (v4) at (.5,0) {};
                \draw (v1) edge +(-.5,0)
                (v2) edge (v3) edge (v1) edge (v4)
                (v3) edge (v1) edge (v4) 
                ;
                \node[int] (w1) at (1.2,0) {};
                \node[int] (w2) at (1.7,.5) {};
                \node[int] (w3) at (1.7,-.5) {};
                \draw (w1) edge (w2) edge (w3) edge (v4)
                (w2) edge (w3) edge +(.5,0)
                (w3) edge +(.5,0);
            \end{tikzpicture}
\mapsto 
\begin{tikzpicture}[scale=.7]
                \node[int] (v1) at (-.5,0) {};
                \node[int] (v2) at (0,.5) {};
                \node[int] (v3) at (0,-.5) {};
                \node[int] (v4) at (.5,0) {};
                \draw (v1) edge +(-.5,0)
                (v2) edge (v3) edge (v1) edge (v4)
                (v3) edge (v1) edge (v4) 
                (v4) edge +(.5,0);
\end{tikzpicture}
\otimes 
\begin{tikzpicture}[scale=.7]
                \node[int] (w1) at (1.2,0) {};
                \node[int] (w2) at (1.7,.5) {};
                \node[int] (w3) at (1.7,-.5) {};
                \draw (w1) edge (w2) edge (w3) edge +(-.5,0)
                (w2) edge (w3) edge +(.5,0)
                (w3) edge +(.5,0);
            \end{tikzpicture}
\end{equation}
Let us define $\Fe(\POp)$ to be the sum over all $g$,
\[
\Fe(\POp)(r) := \bigoplus_{g} \Feyn(\POp)(g,r).
\]
It is a 1-shifted cyclic cooperad, with an additional grading by the genus $g$.
However, one observes that $\Fe(\POp)$ also carries the structure of a cyclic operad. The cyclic operadic composition
\begin{gather*}
   \circ_{i,j}: \Fe(\POp)(r)\otimes \Fe(\POp)(s)
   \to 
   \Fe(\POp)(r+s-2)
\end{gather*}
is obtained by fusing the vertices $u$ (respectively $v$) adjacent to the legs $i$ (respectively $j$).
\[
\circ_{i,j}:
\begin{tikzpicture}
    \node[int, label=90:{$u$}, label=-90:{$\scriptstyle p_u$}] (u) at (0,0) {};
    \node (i) at (.5,0) {$\scriptstyle i$};
    \draw (u) edge (i) edge +(-.5,0) edge +(-.5,0.5) edge +(-.5,-0.5);
\end{tikzpicture}
\otimes 
\begin{tikzpicture}
    \node[int, label=90:{$v$}, label=-90:{$\scriptstyle p_v$}] (u) at (0,0) {};
    \node (j) at (-.5,0) {$\scriptstyle j$};
    \draw (u) edge (j) edge +(.5,0) edge +(.5,0.5) edge +(.5,-0.5);
\end{tikzpicture}
\mapsto 
\begin{tikzpicture}
    \node[int, label=-90:{$\scriptstyle q$}] (u) at (0,0) {};
    \draw (u)  edge +(-.5,0) edge +(-.5,0.5) edge +(-.5,-0.5);
    \draw (u)  edge +(.5,0) edge +(.5,0.5) edge +(.5,-0.5);
\end{tikzpicture}
\]
Here the decoration on the fused vertex is defined as the composition $q := p_u \circ_{i,j} p_v$ of the decorations $p_u$ and $p_v$ on $u$ and $v$. 
It is furthermore obvious that this composition operation cannot create bridges. Hence the cyclic sub-collection 
\[
\Fe(\POp)^{bl} \subset \Fe(\POp).
\]
spanned by bridgeless graphs is a cyclic sub-operad.
Furthermore, any graph splits as a tree of its bridgeless components, and thus we may identify $\Fe(\POp)$ with the cyclic bar construction of the cyclic operad $\Fe(\POp)^{bl}$,
\begin{equation}\label{equ:Fe P B}
\Fe(\POp) \cong B_{cyc} \Fe(\POp)^{bl}.
\end{equation}
With these observations we have already essentially proven our main Theorem.

\begin{thm}\label{thm:main}
Let $\POp$ be a 1-shifted cyclic operad.
Then there is a quasi-isomorphism of cyclic operads
\[
B_{\kk cyc}^c \Fe(\POp) \xrightarrow{\sim} \Fe(\POp)^{bl}.
\]
It gives rise to a spectral sequence 
\[
E_2 = H(B_{\kk cyc}^c H(\Fe(\POp))) \Rightarrow H(\Fe(\POp)^{bl})
\]
that converges to the cohomology of $\Fe(\POp)^{bl}$ as long as $\POp(r)$ is finite dimensional for each $r$ and zero for $r<3$.
\end{thm}

%We may then use that the counit of the cyclic operadic bar-cobar adjunction is a weak equivalence to conclude that we have a quasi-isomorphism
%\[
%B_{cyc}^c \Fe(\POp) \xrightarrow{\sim} \Fe(\POp)^{bl},
%\]
%where $B_{cyc}^c$ is the cyclic operadic cobar construction. From the filtration on the number of vertices we furthermore obtain a spectral sequence
%\[
%E_2 = H(B_{cyc}^c H(\Fe(\POp))) \Rightarrow H(\Fe(\POp)^{bl})
%\]
%converging to the cohomology of the bridgeless complex $\Fe(\POp)^{bl}$. 
%We may also replace cyclic operads by 1-shifted cyclic operads, which is a minor variation defined so that the composition has cohomological degree 1. We may use the straightforward variant of the Feynman transform $\Feyn_{\kk}$ and the cyclic bar construction $B_{\kk cyc}$ for these 1-shifted objects.
%By repeating the above resaoning we then obtain that for $\QOp$ a 1-shifted cyclic operad one has an isomorphism
%\[
%\Fe_{\kk}(\QOp) \cong B_{\kk cyc} \Fe_{\kk}(\QOp)^{bl},
%\]
%and hence also a spectral sequence 
%\[
%E_2 = H(B_{\kk cyc}^c H(\Fe_{\kk}(\QOp))) \Rightarrow H(\Fe_{\kk}(\QOp)^{bl}).
%\]

Now consider the commutative cyclic operad $\Com$ defined such that $\Com(r)=\Q$ for all $r\geq 3$. The Kontsevich graph complex $\G_0$ is the $r=0$ part of the Feynman transform of the commutative operad, and we have
\begin{align*}
\G_0^{g\lp} &\cong \Feyn(\Com)(g,0)
&
\G_0 &\cong \Fe(\Com)(0)
&
\G_0^{bl} &\cong \Fe(\Com)^{bl}(0).
\end{align*}
From Theorem \ref{thm:main} we then immediately obtain:

\begin{cor}
There a a quasi-isomorphism 
    \begin{align*}
        (B_{\kk cyc}^c\Fe(\Com))(0) \to \G_0^{bl}
    \end{align*}
    and a spectral sequence
\[
E_2 = H(B_{\kk cyc}^c H(\Fe(\Com))(0) \Rightarrow H(\G_0^{bl}).
\]
%
%    \begin{align*}
%E_2 = H(B_{cyc}^c H(\Fe(\Com))(0) &\Rightarrow H(\G_0^{bl})
%\\
%    E_2 = H(B_{\kk cyc}^c H(\Fe_{\kk}(\Com_{\kk}))(0)[1] &\Rightarrow H(\G_1^{bl})
%    \end{align*}
    converging to the cohomology of the bridgeless Kontsevich graph complex.
\end{cor}

The above result allows to express the cohomology $H(\G_0^{bl})$ through knowledge on the ordinary Kontsevich graph complex, albeit in its version $\Fe(\Com)$ with external legs. This latter version has been more widely used and studied compared to the bridgeless complex, although we do not fully know its cohomology either.
A slight variant of the above argument also yields parallel results for the complex $\G_1^{bl}$, see section \ref{sec:odd} below.
Finally, the graph complexes $\G_n$ for $n$ of the same parity are isomorphic up to degree shifts
\begin{align}
\label{equ:G inc g}
\G_n^{g\lp} &\cong \G_{n+2}^{g\lp}[-2g]
&
\G_n^{bl,g\lp} &\cong \G_{n+2}^{bl,g\lp}[-2g],
\end{align}
so that all $\G_n^{bl}$ are covered by our results. 
%Numerically, the spectral sequences thus obtained converge at the $E_2$ page up to loop order 10, though we cannot guarantee this in general, see section \ref{sec:numerics} below.

\subsection*{Acknowledgements}
The author is grateful for discussions with Benjamin Brück and Peter Patzt, from which this project arose. 
This work has been partially supported by the NCCR Swissmap funded by the Swiss National Science Foundation, and the Horizon Europe Framework Program CaLIGOLA (101086123).

\section{Background: Operads and the Feynman transform}
We recall here briefly generalities on cyclic operads and the Feynman transform, see \cite{GK} for a more detailed discussion.
A symmetric sequence is a collection $\{\POp(r)\}_{r\geq 0}$ of dg vector spaces so that $\POp(r)$ is equipped with an action of the symmetric group $S_r$. It is also customary to identify symmetric sequences with functors 
\[
\mathcal{F}\mathit{in} \to \mathit{dg}\mathcal{V}\mathit{ect}
\]
from the category of finite sets to dg vector spaces, see \cite[section 1.4]{GK}.
A (non-unital) cyclic operad is a symmetric sequence $\POp$ together with composition operations
\[
\circ_{i,j}\colon \POp(r+1)\otimes \POp(s+1) \to \POp(r+s),
\]
satisfying natural equivariance and associativity axioms.
To a cyclic operad $\POp$ we may associate its Feynman transform $\Feyn(\POp)$, the collection of dg vector spaces $\Feyn(\POp)(g,n)$ as in \eqref{equ:feyn def}. The differential on $\Feyn(\POp)$ has the form $d_{\POp}+d_c$, with $d_{\POp}$ induced from the differential on $\POp$, and $d_c$ is the operation of contracting one edge.
We again refer the reader to \cite{GK} for detailed definitions. We note however that we use slightly different conventions: The Feynman transform of \cite{GK} is the dual of ours, $(\Feyn\POp)^*$.

The collection of dg vector spaces $\Feyn(\POp)$ carries the structure of a dg $\kk$-modular cooperad, i.e., a version of modular cooperad for which the cooperadic coocmposition has cohomological degree +1, see \cite[section 4]{GK}. More precisely, the defining modular cooperadic cocomposition morphisms 
\begin{align*}
\Feyn(\POp)(g,r)[1] & \to \Feyn(\POp)(g',r'+1) \otimes  \Feyn(\POp)(g-g',r-r'+1)
\\
\Feyn(\POp)(g,r)[1] & \to \Feyn(\POp)(g-1,r+2)
\end{align*}
are both given by removing one edge. If the edge is a bridge, the removal gives rise to a cocomposition of the first type (see \eqref{equ:cocomp ex}) and if it is not a bridge, one obtains a cocomposition of the second type.

Every cyclic operad (resp. cooperad) can be considered as a modular operad (resp. cooperad) concentrated in genus $0$. We define a 1-shifted cyclic cooperad to be a $\kk$-modular cooperad concentrated in genus 0. In other words, a 1-shifted cyclic cooperad is the same data as a cyclic cooperad except that the cocomposition has degree +1 and the signs in the axioms are suitably adapted to reflect this change.

There is a variation $\Feyn_{\kk}(\POp)$ of the Feynman transform defined for 1-shifted cyclic operads $\POp$. Similarly to $\eqref{equ:feyn def}$ we have 
\begin{equation}\label{equ:feyn def kk}
\Feyn_{\kk}(\POp)(g,r)= \left(\bigoplus_{\Gamma} \bigotimes_{v\in V(\Gamma) } 
\POp(\mathrm{star}_v) \right)/{\sim}.
\end{equation}
Elements of $\Feyn_{\kk}(\POp)$ are linear combinations of graphs whose vertices are decorated by elements of $\POp$. The difference to $\Feyn(-)$ is that the edges carry cohomological degree 0 instead of degree $-1$.

The cyclic bar construction of a cyclic operad (resp. 1-shifted cyclic operad) $\POp$ is the part of the Feynman transform of genus 0,
\begin{align*}
    (B_{cyc}\POp)(r) &:= \Feyn(\POp)(0,r)
    &
    (B_{\kk cyc}\POp)(r) &:= \Feyn_{\kk}(\POp)(0,r).
\end{align*}
By duality one also defines the cobar construction of a cyclic cooperad (resp. 1-shifted cyclic cooperad) $B_{cyc}^c\COp$ (resp. $B_{\kk cyc}^c\COp$).
One has natural quasi-isomorphisms 
\begin{align}
\label{equ:bar cobar qiso}
B_{\kk cyc}^cB_{cyc} \POp &\xrightarrow{\sim} \POp 
&
B_{cyc}^cB_{\kk,cyc} \POp &\xrightarrow{\sim} \POp.
\end{align}

The Kontsevich graph complex $\G_n$ is the $r=0$-part of the Feynman transform of the commutative cyclic operad $\Com$,
\[
\G_0 := \bigoplus_{g} \Feyn(\Com)(g,0).
\]
It comes with a grading by genus (loop order) $\G_0=\bigoplus_g \G_0^{g\lp}$.
Similarly, we may consider the 1-shifted cyclic cooperad $\Com_{\kk}$ such that
\[
\Com_{\kk}(r) = \begin{cases}
    \Q[3-r] & \text{for $r\geq 3$} \\
    0 & \text{otherwise}
\end{cases}.
\]
One than defines the "odd" Kontsevich graph complex 
\[
\G_{3} := \bigoplus_g \Feyn_{\kk}(\Com_{\kk})(g,0)[-3].
\]
The degree shift by 3 here is merely a convention.
Then finally $\G_n$ may be defined for any integer $n$ via \eqref{equ:G inc g}.

\section{The main results}
\subsection{Proof of Theorem \ref{thm:main}}
We next finish the proof of Theorem \ref{thm:main}. 
Applying the cobar construction $B_{\kk cyc}^c(-)$ to \eqref{equ:Fe P B} and using the quasi-isomorphism \eqref{equ:bar cobar qiso} we obtain the quasi-isomorphism of Theorem \ref{thm:main}
\[
B_{\kk cyc}^c \Fe(\POp) \cong B_{\kk cyc} B_{cyc} \Fe(\POp)^{bl} 
\xrightarrow{\sim} \Fe(\POp)^{bl}.
\]
Note that $B_{\kk cyc}^c$ inherits the grading by loop order from $\Fe(\POp)$. Furthermore, the component of fixed loop order is finite dimensional, by the assumptions of the Theorem that $\POp(r)$ is finite dimensional and zero if $r<3$. 
Next we may equip $B_{\kk cyc}^c \Fe(\POp)$ with another filtration by the number of vertices in trees appearing in the cobar construction and consider the associated spectral sequence. 
This spectral sequence automatically converges to the cohomology because of the finite dimensionality of every genus component.
The differential on $B_{\kk cyc}^c \Fe(\POp)$ is split into two parts
\[
d=d_s + d_{\Fe(\POp)},
\]
with $d_{\Fe(\POp)}$ the differential on $\Fe(\POp)$ and $d_s$ the part of the diffferential of the cobar construction that splits vertices, applying the cocomposition to the decoration in $\Fe(\POp)$. This latter part increases the number of tree-vertices.
The associated graded complex is hence
\[
E_0 \cong (B_{\kk cyc}^c \Fe(\POp), d_{\Fe(\POp)}),
\]
with cohomology 
\[
E_1 = B_{\kk cyc}^c H(\Fe(\POp)).
\]
The differential (induced by $d_s$) on the $E_1$ page is identified with the cobar differential for $B_{\kk cyc}^c H(\Fe(\POp))$. This shows Theorem \ref{thm:main}.\hfill\qed

\subsection{Variant for odd $n$}\label{sec:odd}
The argument of the previous section immediately extends mutatis mutandis to the case of 1-shifted cyclic operads. Thies yields the following results:
\begin{thm}\label{thm:main odd}
Let $\POp$ be a 1-shifted cyclic operad.
There is a quasi-isomorphism of cyclic operads
\[
B_{cyc}^c \Fe_{\kk}(\POp) \xrightarrow{\sim} \Fe_{\kk}(\POp)^{bl}.
\]
It gives rise to a spectral sequence 
\[
E_2 = H(B_{cyc}^c H(\Fe_{\kk}(\POp))) \Rightarrow H(\Fe_{\kk}(\POp)^{bl})
\]
that converges to the cohomology of $\Fe_{\kk}(\POp)^{bl}$ as long as $\POp(r)$ is finite dimensional for each $r$ and zero for $r<3$.
\end{thm}
\begin{cor}
There a a quasi-isomorphism 
    \begin{align*}
        (B_{cyc}^c\Fe_{\kk}(\Com_{\kk}))(0) \to \G_3^{bl}[3]
    \end{align*}
    and a spectral sequence
\[
E_2 = H(B_{cyc}^c H(\Fe_{\kk}(\Com_{\kk}))(0) \Rightarrow H(\G_3^{bl})[3].
\]
%
%    \begin{align*}
%E_2 = H(B_{cyc}^c H(\Fe(\Com))(0) &\Rightarrow H(\G_0^{bl})
%\\
%    E_2 = H(B_{\kk cyc}^c H(\Fe_{\kk}(\Com_{\kk}))(0)[1] &\Rightarrow H(\G_1^{bl})
%    \end{align*}
    converging to the cohomology of the bridgeless Kontsevich graph complex.
\end{cor}

%, and its 1-shifted variant $\Com_{\kk}$ defined such that $\Com((r))=\Q[1-r]$ for $r\geq 3$.
%The Kontsevich graph complex is then just the $n=0$-part of the Feynman transform of these commutative operads,
%\begin{align*}
%\G_0^{g\lp} &\cong \Feyn(\Com)(g,0) \\ 
%\G_1^{g\lp}[-1] &\cong \Feyn_{\kk}(\Com_{\kk})(g,0).
%\end{align*}
%We have hence shown our main result:
%\begin{thm}
%    There are spectral seqences 
%    \begin{align*}
%E_2 = H(B_{cyc}^c H(\Fe(\Com))(0) &\Rightarrow H(\G_0^{bl})
%\\
%    E_2 = H(B_{\kk cyc}^c H(\Fe_{\kk}(\Com_{\kk}))(0)[1] &\Rightarrow H(\G_1^{bl})
%    \end{align*}
%    converging to the bridgeless part of the Kontsevich graph cohomology.
%\end{thm}

\subsection{Variant for simple graphs}
The graph complexes $\G_n$ and $\G_n^{bl}$ above contain non-simple graphs, i.e., graphs with self-loops and/or multiple edges. Denote by $\G_n^s$ (resp. $\G_n^{s,bl}$) the complexes obtained by quotienting $\G_n$ (resp. $\G_n^{bl}$) by the dg ideal spanned by non-simple graphs.
It is well-known \cite{grt,WZ} that the projection $\G_n\to\G_n^s$ is almost a quasi-isomorphism, namely
\begin{align*}
H(\G_0) &\cong H(\G_0^s) \\
H(\G_3) &\cong H(\G_3^s) \oplus \Q \, \thetagr.
\end{align*}

Generally, we may take the quotient of the Feynman transform for any cyclic operad or 1-shifted cyclic operad by the non-simple graphs. Denote the resulting functors $\Feyn^s$, $\Feyn_{\kk}^s$. Applied to the commutative operads this again yields complexes that are known to have almost the same cohomology as their simple counterparts.
For the "even" graph complexes one can extract from the proof of \cite[Lemma 5]{AWZ} that

\begin{align*}
H(\Feyn(\Com)(g,r)) &\cong 
\begin{cases}
    H(\Feyn^s(\Com)(g,r)) & \text{for $(g,r)\neq (1,1)$} \\
    \Q \tadpi \text{for $(g,r)= (1,1)$}
\end{cases}
\\
H(\Feyn^s(\Com)(1,1)) &=0.
\end{align*}

For the "odd" graph complex we have: 
\begin{prop}\label{prop:multi edge}
The projection $\Feyn_{\kk}(\Com_{\kk})(g,r)\to \Feyn_{\kk}^s(\Com_{\kk})(g,r)$ is a quasi-isomorphism for $(g,r)\neq (2,0), (2,1),(1,2)$.
In addition we have 
\begin{align*}
H(\Feyn_{\kk}(\Com_{\kk})(2,0)) &\cong  \Q \, \thetagr \\
H(\Feyn_{\kk}(\Com_{\kk})(2,1)) &\cong  \Q \, \thetagi \\
H(\Feyn_{\kk}(\Com_{\kk})(1,2)) &\cong  \Q \, \thetagii \\
H(\Feyn_{\kk}^s(\Com_{\kk})(2,0)) &=H(\Feyn_{\kk}^s(\Com_{\kk})(2,1))=H(\Feyn_{\kk}^s(\Com_{\kk})(1,2))=0.
\end{align*}
\end{prop}
The proposition can be proven along the lines of \cite[Theorem 2]{WZ}.
However, we provide an independent proof in Appendix \ref{app:multiple}.
Note that in either case, imposing the simplicity condition only alters the cohomology of the Kontsevich graph complex for a finite set of $g,r$, and by some finite dimensional vector space.

Now, repeating the discussion of the previous section for the simple-graph version of the Feynman transform, we obtain:

\begin{thm}\label{thm:main_s}
    There are quasi-isomorphisms
\begin{align*}
B^c_{cyc}(\Feyn^s\Com) (0) &\xrightarrow{\sim} \GC_0^{s,bl} \\
B^c_{\kk cyc}(\Feyn^s_{\kk}\Com_{\kk}) (0) &\xrightarrow{\sim} \GC_3^{s,bl}[3]
\end{align*}
that give rise to spectral sequences
\begin{align*}
E_2 = H(B^c_{cyc}(H(\Feyn^s\Com)) (0)) &\Rightarrow H(\GC_0^{s,bl}) \\
E_2 = H(B^c_{\kk cyc}(H(\Feyn^s_{\kk}\Com_{\kk})) (0)) &\Rightarrow  H(\GC_3^{s,bl}[3]).
\end{align*}
\end{thm}
Note that now the cohomology of $\GC_n^{bl}$ may differ more significantly from its simple variant $\GC_n^{s,bl}$, since non-simple components might appear as constituents of larger trees in the cyclic cobar construction.

\section{Numerical discussion}\label{sec:numerics}
We have implemented the simple version of the bridgeless graph complex on the computer using the GH framework \cite{BW}. A table of the computed cohomology can be found in Figure \ref{fig:cohom}.
For illustration, let us discuss the non-zero classes that appear in the table for $H(\GC_0^{s,bl})$ up to loop order 10.
We recall from the following known facts facts:
\begin{itemize}
    \item By \cite[Proposition 11]{BW} we have that for $g\leq 7$ 
    \[
    H(\Feyn(\Com)(g,1)) \cong H(\Feyn(\Com)(g,0)) =: H(\G_0^{g\lp}).
    \]
    \item By \cite[Figure 7]{BW} we have that 
    \begin{align}\label{equ:g2h2}
        H(\Feyn(\Com)(g,2))
        = 
        \begin{cases}
            \Q \begin{tikzpicture}[scale=.7]
                \node[int] (v1) at (-.5,0) {};
                \node[int] (v2) at (0,.5) {};
                \node[int] (v3) at (0,-.5) {};
                \node[int] (v4) at (.5,0) {};
                \draw (v1) edge +(-.5,0)
                (v2) edge (v3) edge (v1) edge (v4)
                (v3) edge (v1) edge (v4) 
                (v4) edge +(.5,0);
            \end{tikzpicture}
            & \text{for $g=2$} \\
            0 & \text{for $g\neq 2$, $g\leq 4$}
        \end{cases}
    \end{align}
    and 
    \begin{align*}
        H(\Feyn(\Com)(0,3)) &= \Q \, 
            \begin{tikzpicture}
                \node[int] (v) at (0,0) {};
                \draw (v) edge +(-30:.3) edge +(90:.3) edge +(-150:.3);
            \end{tikzpicture}\\
        H(\Feyn(\Com)(1,3)) &= \Q \, 
            \begin{tikzpicture}
                \node[int] (v) at (90:.3) {};
                \node[int] (w) at (-30:.3) {};
                \node[int] (u) at (-150:.3) {};
                \draw (w) edge +(-30:.3) edge (v) edge (u) 
                (v) edge +(90:.3) edge (u)
                (u) edge +(-150:.3);
            \end{tikzpicture}
    \end{align*}
\end{itemize}

Next, let us write down a basis of $B^c_{\kk cyc}(H(\Fe^s(\Com)))(0)$ up to loop order $10$. Since the 1-leg graph cohomology is zero in loop orders $< 3$ it follows that the trees appearing from the cyclic cobar construction can have at most 3 leaves (in loop orders $\leq 10$).
Furthermore, each bivalent vertex in such a tree contributes at least 2 to the genus, so that in the case of two leaves we can have at most 2 bivalent vertices. This yields the following types of trees:
\begin{itemize}
    \item Type $A=\begin{tikzpicture}\node[ext] (v) at (0,0){};\end{tikzpicture}$, i.e., a single vertex. This produces one copy of $H(\GC_0)$ inside $H(\GC_0^{s,bl})$, see \cite[Figure 1]{BW} for a table of $H(\GC_0)$. 
    \item Type $B=\begin{tikzpicture}
    \node[ext] (v) at (0,0){};
    \node[ext] (w) at (0.5,0){};
    \draw (v) edge (w);
    \end{tikzpicture}$. In loop orders $\leq 10$ this produces the second symmetric power of $H(\GC_0)$.
    Concretely, the representatives in $\GC_0^{s,bl}$ of such classes are obtained by joining the two representatives of classes in $H(\GC_0)$ at one vertex, for example:
    
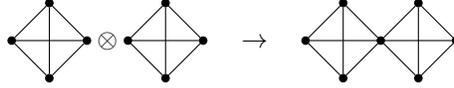
\[
    \begin{tikzpicture}
        \node[int] (v1) at (0,0) {};
        \node[int] (v2) at (1,0) {};
        \node[int] (v4) at (.5,.5) {};
        \node[int] (v5) at (.5,-.5) {};
        \draw (v1) edge (v2) edge (v4) edge (v5) 
        (v4) edge (v5)
        (v2) edge (v4) edge (v5);
    \end{tikzpicture}
    \otimes 
        \begin{tikzpicture}
        \node[int] (v1) at (0,0) {};
        \node[int] (v2) at (1,0) {};
        \node[int] (v4) at (.5,.5) {};
        \node[int] (v5) at (.5,-.5) {};
        \draw (v1) edge (v2) edge (v4) edge (v5) 
        (v4) edge (v5)
        (v2) edge (v4) edge (v5);
    \end{tikzpicture}
    \quad
    \to 
    \quad
\begin{tikzpicture}
        \node[int] (v1) at (0,0) {};
        \node[int] (v2) at (1,0) {};
        \node[int] (v3) at (-1,0) {};
        \node[int] (v4) at (.5,.5) {};
        \node[int] (v5) at (.5,-.5) {};
        \node[int] (v6) at (-.5,.5) {};
        \node[int] (v7) at (-.5,-.5) {};
        \draw (v1) edge (v2) edge (v3) edge (v4) edge (v5) edge (v6) edge (v7)
        (v6) edge (v7) (v4) edge (v5)
        (v2) edge (v4) edge (v5)
        (v3) edge (v6) edge (v7);
    \end{tikzpicture}
    \]
    \item Type $C=\begin{tikzpicture}
    \node[ext] (v) at (0,0){};
    \node[ext] (w) at (0.5,0){};
    \node[ext] (u) at (1,0){};
    \draw (w) edge (v) edge (u);
    \end{tikzpicture}$.
    In loop orders $\leq 10$ this also produces the second symmetric power of $H(\GC_0)$, but shifted in loop order by 2 and in degree by 5, due to the decoration \eqref{equ:g2h2} at the central vertex. The simplest example is
    
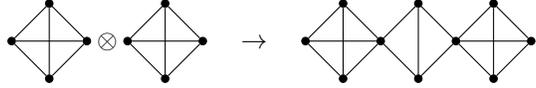
\[
    \begin{tikzpicture}
        \node[int] (v1) at (0,0) {};
        \node[int] (v2) at (1,0) {};
        \node[int] (v4) at (.5,.5) {};
        \node[int] (v5) at (.5,-.5) {};
        \draw (v1) edge (v2) edge (v4) edge (v5) 
        (v4) edge (v5)
        (v2) edge (v4) edge (v5);
    \end{tikzpicture}
    \otimes 
        \begin{tikzpicture}
        \node[int] (v1) at (0,0) {};
        \node[int] (v2) at (1,0) {};
        \node[int] (v4) at (.5,.5) {};
        \node[int] (v5) at (.5,-.5) {};
        \draw (v1) edge (v2) edge (v4) edge (v5) 
        (v4) edge (v5)
        (v2) edge (v4) edge (v5);
    \end{tikzpicture}
    \quad
    \to 
    \quad
\begin{tikzpicture}
        \node[int] (v1) at (0,0) {};
        \node[int] (v2) at (1,0) {};
        \node[int] (v3) at (-1,0) {};
        \node[int] (v4) at (.5,.5) {};
        \node[int] (v5) at (.5,-.5) {};
        \node[int] (v6) at (-.5,.5) {};
        \node[int] (v7) at (-.5,-.5) {};
        \node[int] (w1) at (2,0) {};
        \node[int] (w2) at (1.5,0.5) {};
        \node[int] (w3) at (1.5,-.5) {};
        \draw (v1) edge (v3) edge (v4) edge (v5) edge (v6) edge (v7)
        (v6) edge (v7) (v4) edge (v5)
        (v2) edge (v4) edge (v5)
        (v3) edge (v6) edge (v7)
        (w1) edge (w2) edge (w3) edge (v2)
        (w2) edge (w3) edge (v2) 
        (w3) edge (v2);
    \end{tikzpicture}.
    \]
    \item Type $D=\begin{tikzpicture}
    \node[ext] (v) at (0,0){};
    \node[ext] (w) at (0.5,0){};
    \node[ext] (u) at (-.5,0){};
    \node[ext] (x) at (0,0.5){};
    \draw (v) edge (x) edge (w) edge (u);
    \end{tikzpicture}$. Here the only possible decorations are the three-loop classes (in even degree) at the leaves, and the 0-loop class at the central vertex.
    
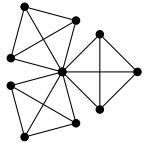
\[
    \begin{tikzpicture}%[transform shape,rotate around={0:(0,0)}]
    \begin{scope}[rotate around={120:(0,0)}]
        \node[int] (v1) at (0,0) {};
        \node[int] (v2) at (1,0) {};
        \node[int] (v4) at (.5,.5) {};
        \node[int] (v5) at (.5,-.5) {};
        \draw (v1) edge (v2) edge (v4) edge (v5) 
        (v4) edge (v5)
        (v2) edge (v4) edge (v5);
    \end{scope}
    \begin{scope}[rotate around={-120:(0,0)}]
        \node[int] (v1) at (0,0) {};
        \node[int] (v2) at (1,0) {};
        \node[int] (v4) at (.5,.5) {};
        \node[int] (v5) at (.5,-.5) {};
        \draw (v1) edge (v2) edge (v4) edge (v5) 
        (v4) edge (v5)
        (v2) edge (v4) edge (v5);
    \end{scope}
        \begin{scope}[]
        \node[int] (v1) at (0,0) {};
        \node[int] (v2) at (1,0) {};
        \node[int] (v4) at (.5,.5) {};
        \node[int] (v5) at (.5,-.5) {};
        \draw (v1) edge (v2) edge (v4) edge (v5) 
        (v4) edge (v5)
        (v2) edge (v4) edge (v5);
    \end{scope}
    \end{tikzpicture}
    \]
    \item Type $E=\begin{tikzpicture}
    \node[ext] (v) at (0,0){};
    \node[ext] (w) at (0.5,0){};
    \node[ext] (u) at (1,0){};
    \node[ext] (x) at (1.5,0){};
    \draw (u) edge (x) edge (w) (v) edge (w);
    \end{tikzpicture}$. Here the only possible decoration are the three-loop classes (in even degree) at the leaves, and the 2-loop classes (in odd degree) at the bivalent vertices. Hence the graph vanishes by symmetry and there is no contribution.
    
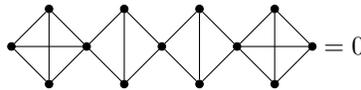
\[
    \begin{tikzpicture}
        \node[int] (v1) at (0,0) {};
        \node[int] (v2) at (1,0) {};
        \node[int] (v3) at (-1,0) {};
        \node[int] (v4) at (.5,.5) {};
        \node[int] (v5) at (.5,-.5) {};
        \node[int] (v6) at (-.5,.5) {};
        \node[int] (v7) at (-.5,-.5) {};
        \node[int] (w1) at (2,0) {};
        \node[int] (w2) at (1.5,0.5) {};
        \node[int] (w3) at (1.5,-.5) {};
        \node[int] (u1) at (3,0) {};
        \node[int] (u2) at (2.5,0.5) {};
        \node[int] (u3) at (2.5,-.5) {};        
        \draw (v1) edge (v3) edge (v4) edge (v5) edge (v6) edge (v7)
        (v6) edge (v7) (v4) edge (v5)
        (v2) edge (v4) edge (v5)
        (v3) edge (v6) edge (v7)
        (w1) edge (w2) edge (w3) 
        (w2) edge (w3) edge (v2) 
        (w3) edge (v2)
        (u1) edge (u2) edge (u3) edge (w1)
        (u2) edge (u3) edge (w1) 
        (u3) edge (w1);
    \end{tikzpicture}
    =0
    \]
\end{itemize}

Collecting all such basis elements one sees that they account for all of $H(\GC_0^{s,bl})$ in loop orders $\leq 10$, and there cannot be further cancellations, so that our spectral sequence abuts at $E_1$ in loop orders $\leq 10$.
The classes we found above are listed in Figure \ref{fig:cohom 2} for illustration. There an entry "$3=1_A+2_B$" means that of the 3-dimensional cohomology, one dimension is spanned by classes of type $A$, and two by classes of type $B$.
Also note that the table lists classes by vertex number, with the vertex number $v$ determined from the degree $d$ and the loop order $g$ by the formula
\[
v = -d-g+1.
\]
Furthermore, the degree is just the sum of the decorations of the vertices of the tree appearing in the cobar construction -- the edges in that cobar construction carry degree 0. In turn, the degree of those decorations (i.e., graphs) is minus the number of edges in the graphs. 

\begin{figure}

%\scalebox{ 1 }{
 \begin{tabular}{|g|D|D|D|D|D|D|D|D|D|D|D|D|D|D|D|D|D|}
  \hline
  \rowcolor{Gray}
    $g$,$v$ & 4 & 5 & 6 & 7 & 8 & 9 & 10 & 11 & 12 & 13 & 14 & 15 & 16 & 17 & 18 & 19 \\ \hline
3 &  1 &  - &  - &  - &  - &  - &  - &  - &  - &  - &  - &  - &  - &  - &  - &  - \\
\hline
4 &  - &  0 &  1 &  - &  - &  - &  - &  - &  - &  - &  - &  - &  - &  - &  - &  - \\
\hline
5 &  - &  0 &  0 &  0 &  2 &  - &  - &  - &  - &  - &  - &  - &  - &  - &  - &  - \\
\hline
6 &  - &  0 &  0 &  1 &  0 &  0 &  2 &  - &  - &  - &  - &  - &  - &  - &  - &  - \\
\hline
7 &  - &  - &  0 &  0 &  0 &  2 &  0 &  0 &  3 &  - &  - &  - &  - &  - &  - &  - \\
\hline
8 &  - &  - &  0 &  0 &  0 &  0 &  1 &  4 &  0 &  0 &  4 &  - &  - &  - &  - &  - \\
\hline
9 &  - &  - &  0 &  0 &  0 &  0 &  1 &  0 &  2 &  7 &  0 &  0 &  5 &  - &  - &  - \\
\hline
10 &  - &  - &  0 &  0 &  0 &  0 &  0 &  0 &  2 &  0 &  5 &  11 &  0 &  0 &  6 &  - \\
\hline
\end{tabular}

\medskip

 \begin{tabular}{|g|D|D|D|D|D|D|D|D|D|D|D|D|D|D|D|D|D|}
  \hline
 \rowcolor{Gray}
$g$,$v$ & 4 & 5 & 6 & 7 & 8 & 9 & 10 & 11 & 12 & 13 & 14 & 15 & 16 & 17 & 18 & 19 \\ \hline
3 &  1 &  - &  - &  - &  - &  - &  - &  - &  - &  - &  - &  - &  - &  - &  - &  - \\
\hline
4 &  - &  0 &  0 &  - &  - &  - &  - &  - &  - &  - &  - &  - &  - &  - &  - &  - \\
\hline
5 &  - &  0 &  1 &  0 &  0 &  - &  - &  - &  - &  - &  - &  - &  - &  - &  - &  - \\
\hline
6 &  - &  0 &  0 &  1 &  0 &  0 &  1 &  - &  - &  - &  - &  - &  - &  - &  - &  - \\
\hline
7 &  - &  - &  0 &  0 &  1 &  0 &  0 &  0 &  0 &  - &  - &  - &  - &  - &  - &  - \\
\hline
8 &  - &  - &  0 &  0 &  0 &  2 &  1 &  0 &  1 &  0 &  0 &  - &  - &  - &  - &  - \\
\hline
9 &  - &  - &  0 &  0 &  0 &  0 &  2 &  0 &  0 &  2 &  0 &  0 &  0 &  - &  - &  - \\
\hline
10 &  - &  - &  0 &  0 &  0 &  0 &  0 &  3 &  1 &  0 &  2 &  0 &  0 &  0 &  1 &  - \\
\hline
\end{tabular}
\caption{\label{fig:cohom} Cohomology dimensions of the bridgeless simple graph complex $\GC_n^{s,bl}$ for $n$ odd (top) and $n$ even (bottom). The rows correspond to the loop orders, the columns to the numbers of vertices.
Entries "-" are zero because the complex is zero in those degrees. }
\end{figure}

\begin{figure}
\scalebox{.7}{
 \begin{tabular}{|g|D|D|D|D|D|D|D|D|D|D|D|D|D|D|D|D|}
  \hline
 \rowcolor{Gray}
$g$,$v$ & 4 & 5 & 6 & 7 & 8 & 9 & 10 & 11 & 12 & 13 & 14 & 15 & 16 & 17 & 18  \\ \hline
3 &  $1=1_A$ &  - &  - &  - &  - &  - &  - &  - &  - &  - &  - &  - &  - &  - &  -  \\
\hline
4 &  - &  0 &  0 &  - &  - &  - &  - &  - &  - &  - &  - &  - &  - &  - &  -  \\
\hline
5 &  - &  0 &  $1=1_A$ &  0 &  0 &  - &  - &  - &  - &  - &  - &  - &  - &  - &  -  \\
\hline
6 &  - &  0 &  0 &  $1=1_B$ &  0 &  0 &  $1=1_A$ &  - &  - &  - &  - &  - &  - &  - &  -  \\
\hline
7 &  - &  - &  0 &  0 &  $1=1_A$ &  0 &  0 &  0 &  0 &  - &  - &  - &  - &  - &  -  \\
\hline
8 &  - &  - &  0 &  0 &  0 &  $2=1_A+1_B$ &  $1=1_C$ &  0 &  $1=1_A$ &  0 &  0 &  - &  - &  - &  -  \\
\hline
9 &  - &  - &  0 &  0 &  0 &  0 &  $2=1_A+1_D$ &  0 &  0 &  $2=1_A+1_B$ &  0 &  0 &  0 &  - &  -  \\
\hline
10 &  - &  - &  0 &  0 &  0 &  0 &  0 &  $3=1_A+2_B$ &  $1=1_C$ &  0 &  $2=2_A$ &  0 &  0 &  0 &  $1=1_A$  \\
\hline
\end{tabular}
}
\caption{\label{fig:cohom 2} A copy of the second table of Figure \ref{fig:cohom}, but with the different cohomology classes associated to trees of the cyclic bar construction as explained in the text.}
\end{figure}

\begin{comment}
\begin{figure}
\scalebox{.5}{
 \begin{tabular}{|g|D|D|D|D|D|D|D|D|D|D|D|D|D|D|D|D|D|}
  \hline
  \rowcolor{Gray}
    $g$,$v$ & 4 & 5 & 6 & 7 & 8 & 9 & 10 & 11 & 12 & 13 & 14 & 15 & 16 & 17 & 18 & 19 \\ \hline
3 &  $1=1_A$ &  - &  - &  - &  - &  - &  - &  - &  - &  - &  - &  - &  - &  - &  - &  - \\
\hline
4 &  - &  0 &  $1=1_A$ &  - &  - &  - &  - &  - &  - &  - &  - &  - &  - &  - &  - &  - \\
\hline
5 &  - &  0 &  0 &  0 &  $2=2_A$ &  - &  - &  - &  - &  - &  - &  - &  - &  - &  - &  - \\
\hline
6 &  - &  0 &  0 &  $1=1_A$ &  0 &  0 &  $2=2_A$ &  - &  - &  - &  - &  - &  - &  - &  - &  - \\
\hline
7 &  - &  - &  0 &  0 &  0 &  $2=1_A+1_B$ &  0 &  0 &  $3=3_A$ &  - &  - &  - &  - &  - &  - &  - \\
\hline
8 &  - &  - &  0 &  0 &  0 &  0 &  $1=1_C$ &  $4=2_A+2_B$ &  0 &  0 &  $4=4_A$ &  - &  - &  - &  - &  - \\
\hline
9 &  - &  - &  0 &  0 &  0 &  0 &  $1=1_B$ &  0 &  $2=2_C$ &  $7=3_A+4_B$ &  0 &  0 &  $5=5_A$ &  - &  - &  - \\
\hline
10 &  - &  - &  0 &  0 &  0 &  0 &  0 &  0 &  $2=2_B$ &  0 &  $5=5_C$ &  $11=5_A+6_B$ &  0 &  0 &  $6=6_A$ &  - \\
\hline
\end{tabular}
}
\caption{\label{fig:cohom 3} A copy of the first table of Figure \ref{fig:cohom}, but with the different cohomology classes associated to trees of the cyclic bar construction as explained in the text.}
\end{figure}
\end{comment}

\appendix 

\section{Proof of Proposition \ref{prop:multi edge}}
\label{app:multiple}
Elements of $A_{g,r}:= \Feyn_{\kk}(\Com_{\kk})(g,r)$ are linear combinations of connected graphs with $g$ loops and $r$ numbered external legs. All vertices must have valence $\geq 3$.
These graphs may have multiple edges. However, by symmetry they may not have self-edges, since the isomorphism of the graph that flips the self-edge acts with a sign.
\[
\begin{tikzpicture}
    \node[int] (v) at (0,0) {};
    \draw (v) edge[loop] (v) edge +(-.3,-.3) edge +(0,-.3) edge +(.3,-.3);
\end{tikzpicture}
=
-\begin{tikzpicture}
    \node[int] (v) at (0,0) {};
    \draw (v) edge[loop] (v) edge +(-.3,-.3) edge +(0,-.3) edge +(.3,-.3);
\end{tikzpicture}
=0
\]
The complex $B_{g,r}:= \Feyn_{\kk}^s(\Com_{\kk})(g,r)$ is obtained from $A_{g,r}$ by setting all graphs with multiple edges to zero, and it is our goal to check that the projection $A_{g,r}\to B_{g,r}$ is a quasi-isomorphism for all $(g,r)\neq (2,0), (1,2), (2,1)$. 
To this end we will introduce an auxiliary complex $C_{g,r}$ through which the projection $A_{g,r}\to B_{g,r}$ factors,
\[
A_{g,r} \xrightarrow{\alpha} C_{g,r} \xrightarrow{\beta} B_{g,r},
\]
such that that $\alpha$ and $\beta$ are both quasi-isomorphisms.
Concretely, $C_{g,r}$ is the complex spanned by connected graphs with two kinds of edges, normal edges and fat edges, with $r$ external legs and genus $g$.
Here the genus is the loop order plus the number of fat edges,
\[
g= \#\text{loops} + \#\text{fat edges}.
\]
We require all vertices to be at least trivalent, but we compute the valence of a vertex $v$ as the number of incident normal half-edges, plus \emph{twice} the number of incident fat half-edges.
Finally, we do not allow graphs with fat self-edges, or equivalently set those graphs to zero:
\[
\begin{tikzpicture}
    \node[int] (v) at (0,0) {};
    \draw (v) edge[loop, very thick] (v) edge +(-.3,-.3) edge +(0,-.3) edge +(.3,-.3);
\end{tikzpicture}
=0
\]
The differential on $C_{g,r}$ has the form 
\[
d=d_c+d_{fat},
\]
with $d_c$ summing over all ways of contracting a normal edge as before, and $d_{fat}$ replacing a fat edge by a pair of parallel normal edges 
\[
d_{fat} : 
\begin{tikzpicture}
    \node[int] (v) at (0,0) {};
    \node[int] (w) at (0.7,0) {};
    \draw (v) edge[ very thick] (w) edge +(-.3,-.3) edge +(-.3,0) edge +(-.3,.3)
    (w) edge +(.3,-.3) edge +(.3,0) edge +(.3,.3);
\end{tikzpicture}
\mapsto 
\begin{tikzpicture}
    \node[int] (v) at (0,0) {};
    \node[int] (w) at (0.7,0) {};
    \draw (v) edge[ bend left] (w) edge[ bend right] (w) edge +(-.3,-.3) edge +(-.3,0) edge +(-.3,.3)
    (w) edge +(.3,-.3) edge +(.3,0) edge +(.3,.3);
\end{tikzpicture}
\]
One easily checks that $d^2=0$.
The morphism $\alpha:A_{g,r}\to C_{g,r}$ is just the inclusion.
The morphism $\beta:C_{g,r}\to B_{g,r}$ is the projection obtained by sending all graphs with fat edges or multiple edges to zero, and other graphs to themselves, understood as elements of $B_{g,r}$. It is clear that both $\alpha$, $\beta$ respect the differentials.

\medskip

\underline{Claim 1: $\beta$ is a quasi-isomorphism.}
Consider the filtrations on $B_{g,r}$ and $C_{g,r}$ by the number of vertices of graphs. It is sufficient to show that the morphism induced by $\beta$ on the associated graded complexes 
\[
(C_{g,r}, d_{fat}) \to (B_{g,r}, 0)
\]
is a quasi-isomorphism. 
To this end we define a homotopy $h:C_{g,r}\to C_{g,r}$ by summing over all ways of replacing a $k$-fold normal edge between two vertices (with $k\geq 2$) by a fat edge parallel to $k-2$ normal edges.
\begin{align*}
    h: 
    \begin{tikzpicture}
    \node[int] (v) at (0,0) {};
    \node[int] (w) at (0.7,0) {};
    \node at (.35,.3) {$\scriptstyle k\times$};
    \draw (v) edge (w) edge[ bend left] (w) edge[ bend right] (w) edge +(-.3,-.3) edge +(-.3,0) edge +(-.3,.3)
    (w) edge +(.3,-.3) edge +(.3,0) edge +(.3,.3);
\end{tikzpicture}
\mapsto 
    \begin{tikzpicture}
    \node[int] (v) at (0,0) {};
    \node[int] (w) at (0.7,0) {};
    \node at (.35,.3) {$\scriptstyle k-2\times$};
    \draw (v)  edge[ bend left] (w) edge[very thick, bend right] (w) edge +(-.3,-.3) edge +(-.3,0) edge +(-.3,.3)
    (w) edge +(.3,-.3) edge +(.3,0) edge +(.3,.3);
\end{tikzpicture}.
\end{align*}
Then one easily checks that 
\[
(d_{fat}h+hd_{fat})\gamma = N(\gamma) \gamma,
\]
where $N(\gamma)$ is the number of pairs of vertices of $\gamma$ that are connected by a fat edge or a multiple normal edge.
Hence by \cite[Lemma 7]{FFW} the associated graded of $\beta$ is a quasi-isomorphism, and hence also $\beta$ is a quasi-isomorphism as well.

\underline{Claim 2: $\alpha$ is a quasi-isomorphism as long as $(g,r)\neq (2,0),(2,1),(1,2)$.}
We have to check that the quotient $C_{g,r}/A_{g,r}$ is acyclic. 
Considering the filtration by the number of fat edges, it is in fact sufficient to check that the associated graded complex is acyclic,
\[
H(C_{g,r}/A_{g,r}, d_c)=0. 
\]
Also note that 
\[
(C_{g,r}/A_{g,r}, d_c) \cong \bigoplus_{k\geq 1} (C_{g,r}^{(k)}, d_c),
\]
with $C_{g,r}^{(k)}\subset C_{g,r}$ being the subspace spanned by graphs with exactly $k$ fat edges.
We may temporarily pass to a larger complex $(C_{g,r}^{k}, d_c)$, in which all the $k$ edges are numbered and directed, so that 
\[
C_{g,r}^{(k)} \cong (C_{g,r}^{k})_{S_2\wr S_k}
\]
is identified with the coinvariants under the natural action of the wreath product $S_2\wr S_k$ by renumbering fat edges and changing fat edge directions.
Now we may use a trick of Lambrechts-Volic \cite{LV} and split 
\[
C_{g,r}^{k} \cong 
\begin{tikzcd}[column sep =.5]
U_1 \ar[loop above]{} \ar[bend left]{rr}{d_c'} & \oplus &  \ar[loop above]{} U_2
\end{tikzcd}
\]
with $U_1$ spanned by graphs for which the endpoint $v$ of fat edge 1 has valence $3$, and $U_2$ spanned by all other graphs.
Note that $v$ having valence 3 means that on top of fat edge 1 ending at $v$ (contributing 2 to the valence by our counting) there is exactly one normal edge or hair incident at vertex $v$.
The piece of the differential $d_c':U_1\to U_2$ contracts the unique normal edge at vertex $v$, if any is present. The key point is now that $d_c'$ is surjective, with a one-sided inverse given by splitting off an edge at $v$.
It follows that one has a quasi-isomorphism (see \cite[Lemma 2.1]{PW})
\[
C_{g,r}^{k} \simeq \ker d_c'.
\]
The kernel of $d_c'$ is spanned by graphs in $U_1$ for which there is no normal edge at $v$ that can be contracted. 
This is the case if either $v$ carries an external leg, or $v$ is connected to the other endpoint $w$ of fat edge $1$. 
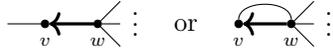
\begin{align*}
    \begin{tikzpicture}
     \node[int,label=-90:{$\scriptstyle v$}] (v) at (0,0) {};   
     \node[int,label=-90:{$\scriptstyle w$}] (w) at (.7,0) {}; 
     \node at (1.2,0.1) {$\vdots$};
     \draw (v) edge[very thick,<-] (w) edge +(-.5,0)
     (w) edge +(.3,.3) edge +(.3,0) edge +(.3,-.3) ;
    \end{tikzpicture}
    \quad \text{or} \quad 
    \begin{tikzpicture}
     \node[int,label=-90:{$\scriptstyle v$}] (v) at (0,0) {};   
     \node[int,label=-90:{$\scriptstyle w$}] (w) at (.7,0) {}; 
     \node at (1.2,0.1) {$\vdots$};
     \draw (v) edge[very thick,<-] (w) edge[out=90, in=90] (w)
     (w) edge +(.3,.3) edge +(.3,0) edge +(.3,-.3) ;
    \end{tikzpicture}
\end{align*}

Next, we may proceed for the complex $\ker d_c'$ in the same manner and split it into subspaces 
\[
C_{g,r}^{k} \cong 
\begin{tikzcd}[column sep =.5]
Y_1 \ar[loop above]{} \ar[bend left]{rr}{d_c''} & \oplus &  \ar[loop above]{} Y_2
\end{tikzcd}
\]
where $Y_1$ is spanned by graphs such that vertex $w$ has exactly one incident normal edge or leg, not counting the edges connecting to $v$, and $Y_2$ is spanned by all remaining graphs.
The part of the differential $d_c'':Y_1\to Y_2$ contracts the unique normal edge incident at $w$ not connecting to $v$, if such an edge is present. As before, $d_c''$ is a surjective map, so that again by \cite[Lemma 2.1]{PW} we have the quasi-isomorphism  
\[
\ker d_c' \simeq \ker d_c''.
\]
Now $\ker d_c''$ is given by graphs with no contractible edge at vertex $w$, and there are only the following three:
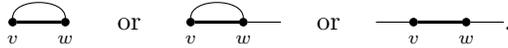
\[
    \begin{tikzpicture}
     \node[int,label=-90:{$\scriptstyle v$}] (v) at (0,0) {};   
     \node[int,label=-90:{$\scriptstyle w$}] (w) at (.7,0) {}; 
     \draw (v) edge[very thick] (w) edge[out=90, in=90] (w)
     ;
    \end{tikzpicture}
    \quad \text{ or } \quad
    \begin{tikzpicture}
     \node[int,label=-90:{$\scriptstyle v$}] (v) at (0,0) {};   
     \node[int,label=-90:{$\scriptstyle w$}] (w) at (.7,0) {}; 
     \draw (v) edge[very thick] (w) edge[out=90, in=90] (w)
     (w) edge +(.5,0) ;
    \end{tikzpicture}
    \quad \text{ or } \quad
    \begin{tikzpicture}
     \node[int,label=-90:{$\scriptstyle v$}] (v) at (0,0) {};   
     \node[int,label=-90:{$\scriptstyle w$}] (w) at (.7,0) {}; 
     \draw (v) edge[very thick] (w) edge +(-.5,0)
     (w) edge +(.5,0);
    \end{tikzpicture}.
\]
These graphs live in $(g,r)=(2,0)$, $(2,1)$ and $(1,2)$. Hence in all other cases we have shown that $\alpha$ is a quasi-isomorphism.

By Claim 1 and Claim 2 above, the projection $A_{g,r}\to B_{g,r}$ is a quasi-isomorphism as long as $(g,r)\neq (2,0),(2,1),(1,2)$, since this projection agrees with the composition $\beta\circ \alpha$.

Finally, the statements of Proposition \ref{prop:multi edge} for $(g,r)= (2,0), (1,2), (2,1)$ are already true on the level of complexes, i.e., 
\begin{align*}
A_{2,0} &\cong  \Q \, \thetagr \\
A_{2,1} &\cong  \Q \, \thetagi \\
A_{1,2} &\cong  \Q \, \thetagii \\
B_{2,0} &= B_{2,1} =B_{1,2}=0.
\end{align*}
The statements on cohomology hence trivially follow.\hfill\qed

\end{document}